\documentstyle{article}
\title{A nonstandard proof of the Jordan curve theorem}
\author{Vladimir Kanovei~\thanks{The author 
acknowledges the support of AMS, DFG, and 
University of Wuppertal.}
\thanks{Moscow Transport Engineering Institute, Russia, \ 
{\tt kanovei@math.uni-wuppertal.de} \ and \ 
{\tt kanovei@mech.math.msu.su}}
\and
Michael Reeken~\thanks 
{Bergische Universit\"at GHS Wuppertal, Germany, \ 
{\tt reeken@math.uni-wuppertal.de}}
}
\date{August 1996}

\font\tenmsb=msbm10 
\font\bfit=cmbxti10
\newcommand{\rbox}[1]{{\rm{#1}}} 
\newcommand{\emph}[1]  {{\em #1}}   
\newcommand{\textsc}[1]{{\sc #1}}       
\newcommand{\textit}[1]{{\it #1}}   
\newcommand{\textbf}[1]{{\bf #1}}   
\newtheorem{jcth}{The Jordan curve theorem}

\newtheorem{lemma}{Lemma}
\newtheorem{proposition}[lemma]{Proposition}
\newtheorem{remark}{\bfit Comment.}
\newtheorem{prooF}{Proof}

\newcommand{\ble} {\begin{lemma} }
\newcommand{\ele} {\end{lemma}}
\newcommand{\bpro}{\begin{proposition} } 
\newcommand{\epro}{\end{proposition}} 
\newcommand{\bpf} {\begin{prooF}\rm } 
\newcommand{\epf} {\qed\end{prooF}} 
\newcommand{\brem}{\begin{remark}\rm }

\newcommand{\erem}{\end{remark}} 
\newcommand{\qed} {\hfill$\msur\Box\msur$} 
\renewcommand{\section}{\subsection}

\newcommand{\ben}{\begin{enumerate}}
\newcommand{\een}{\end{enumerate}}
\newcommand{\IST}{\textbf{IST}}
\newcommand{\al} {\alpha}
\newcommand{\ba} {\beta}
\newcommand{\ga} {\gamma}
\newcommand{\da} {\delta}
\newcommand{\Da} {\Delta}
\newcommand{\vep}{\varepsilon}

\newcommand{\sg} {\sigma}
\newcommand{\cD}{{\cal D}}
\newcommand{\cL}{{\cal L}}
\newcommand{\cR}{{\cal R}}
\newcommand{\cK}{{\cal K}}
\newcommand{\aK}{{{}^\ast\hspace{-0.3ex}\cK}}
\newcommand{\dR}{\hbox{\tenmsb R}}
\newcommand{\<}   {\leq}

\newcommand{\ti}  {\times}
\newcommand{\lra} {\longrightarrow}
\newcommand{\dm}  {$$}
\newcommand{\ans} [1]{\{\hspace{0.2mm}#1\hspace{0.2mm}\}}
\newcommand{\inte}{_{\rbox{int}}}
\newcommand{\exte}{_{\rbox{ext}}}
\newcommand{\new} {_{\rbox{new}}}
\newcommand{\dist}{{\rbox d}}

\newcommand{\ie}  {\hbox{\it i.\hspace*{0.4ex}e\/}.}
\newcommand{\etc} {\hbox{\it etc\/}.}
\newcommand{\dd}[2]{\hpsur\hbox{\mathsurround=0mm${#1}$-#2}}
\newcommand{\noi} {\noindent}
\newcommand{\vom} {\vspace{1mm}}
\newcommand{\wed}[2] {#1\hspace*{0.1ex}#2}
\newcommand{\curve}[2]{\widehat{#1\hspace{0.1ex}
#2\hspace{-2ex}}\hspace{2ex}}%

\newcommand{\msur}{\hspace{-1\mathsurround}}
\newcommand{\hsur}{\hspace{-0.5\mathsurround}}
\newcommand{\hpsur}{\hspace{0.3\mathsurround}}

\newcommand{\mon}{{\rbox{\tt monad}}\,}

\begin{document}
\maketitle  
\begin{abstract}
We give a nonstandard variant of Jordan's proof of the Jordan 
curve theorem which is free of the defects his contemporaries 
criticized and avoids the epsilontic burden of the classical 
proof. The proof is 
self-contained, except for the Jordan theorem for polygons 
taken for granted. 
\end{abstract}

%\hspep
\section*{Introduction}

\noi
The \emph{Jordan curve theorem\/} 
%(abbreviated as JCT)\cite{Jo1893} 
was one of the starting points in the modern development of 
\emph{topology\/} (originally called \emph{Analysis Situs}). 
This result is considered difficult to prove,
at least compared to its intuitive evidence.

\textsc{C.\ Jordan} \cite{Jo1893} considered the assertion to 
be evident for simple polygons and reduced the case of a simple 
closed continuous curve to that of a polygon by approximating the 
curve by a sequence of suitable simple polygons.

Although the idea appears natural to an analyst it is not so easy 
to carry through. \textsc{Jordan}'s proof did not satisfy 
mathematicians of his time. On one hand it was felt that the case 
of polygons also needed a proof based on clearly stated geometrical 
principles, on the other hand his proof was considered incomplete 
(see the criticisms formulated in \cite{Ve05} and in \cite{Os12}).

If one is willing to assume slightly more than mere continuity of 
the curve than much simpler proofs (including the case of polygons) 
are available (see \textsc{Ames} \cite{Am04} and \textsc{Bliss} 
\cite{Bl04} under restrictive hypotheses).

\textsc{O.\ Veblen} \cite{Ve05} is considered the first to have 
given a rigorous proof which, in fact, makes no use
of metrical properties, or, in the words of \textsc{Veblen}: 
\emph{We accordingly assume nothing about analytic geometry, the 
parallel axiom, congruence relations, nor the existence of points 
outside a plane. }

His proof is based on the incidence and order axioms for the plane 
and the natural topology defined by the basis consisting 
of nondegenerate triangles. He also defines simple curves 
intrinsically as specific sets without parametrizations by intervals 
of the real line. He finally discusses how the introduction of one 
additional axiom, existence of a point outside the plane, allows to 
reduce his result to the context \textsc{Jordan} was working in.

\textsc{Veblen} also gave a specific proof for polygons based on 
the incidence and order axioms exclusively (see \cite{Ve04}) which was 
later criticized as inconclusive by \textsc{H.~Hahn} \cite{Ha08} who 
published his own version of a proof based on \textsc{Veblen}'s 
incidence and order axioms of the plane (which, by the way, are 
equivalent to the incidence and order axioms of \textsc{Hilbert}'s 
system).

\textsc{Jordan}'s proof in his \emph{Cours d' analyse }of 1893 is 
elementary as to the tools employed. Nevertheless the proof extends 
over nine pages and, as mentioned above, cannot be considered 
complete. We are interested here in this proof. It depends on some 
facts for polygons and an approximation argument. It is, therefore, 
a natural idea to use nonstandard arguments to eliminate the 
epsilontic burden of the approximation.

There is an article by \textsc{L.\ Narens} \cite{Na71} in which 
this point of view is adopted. Unfortunately, beside some problems 
discovered in his proof 
%(as was communicated by \textsc{W. Henson} 
%on the mailbox {\sc Galaxy}, letter of Oct.\ 31, 1995). 
by those who read is carefully, the reasoning involves more 
complicated topological ideas and anyway is not essentially 
shorter than or comparably elementary as the \textsc{Jordan}'s proof.
%~\footnote
%{\ Cited from the mentioned e-mail of Oct.\ 31, 1995 : 
%\textit{The case which the author\/} (Narens) 
%\textit{fails to take care of seems in fact to be very characteristic 
%of the situation which makes the Jordan Curve Theorem difficult to 
%prove by elementary topological means. ... 
%The basic approach of Narens is appealing -- you approximate the 
%curve by a hyperfinite polygon and take advantage of the fact that 
%the JCT is trivial for polygons. But it seems hard to make this 
%approach work.}}

It is certainly true that not all classical arguments can be replaced 
in some useful or reasonable way by simpler nonstandard arguments. 
But as we shall show it is possible to simplify the approximation 
argument specific to \textsc{Jordan}'s proof. We shall follow the 
proof quite closely but take a
somewhat different approach when proving path-connectedness. 

That nonstandard analysis can even give some additional insight 
into the geometric problem is manifest from the proof by 
\textsc{N.\ Bertoglio} and \textsc{R.\ Chuaqui} \cite{BeCh94} which 
avoids polygons and approximations
entirely by looking at a nonstandard discretization of the plane and
reducing the problem to a combinatorial version of the JCT  
proved by \textsc{L.\ N.\ Stout} \cite{St88}. This reduction of the 
problem to a (formally) discrete one is interesting and leads to a 
proof which establishes a link to a context totally different from 
\textsc{Jordan}'s.

As a curiosity we note in passing that \textsc{Jordan} speaks of 
\emph{infinitesimals\/} in his proof but it is only a figure of 
speech for a number which may be chosen as small as one wishes or 
for a function which tends to zero. 

For reference we state:

\begin{jcth}
\emph{(abbreviated as JCT)}\\[0pt] 
A simple closed continuous curve\/ $\cK$ in the plane separates 
its complement into two open sets of which it is the common 
boundary$;$ one of them is called the \emph{outer region 
(}or \emph{exterior domain)} $\cK\exte$ which is an open, 
unbounded, path-connected set and another set called the 
\emph{inner region (}or \emph{interior domain)} $\cK\inte$ which 
is an open, simply path-connected, bounded set.
\end{jcth}
{\bfit Acknowledgements} \ The authors acknowledge with pleasure 
useful discussions with C.\ W.\ Henson, H.\ J.\ Keisler, S.\ Leth, 
P.\ A.\ Loeb in matters of the Jourdan curve theorem in the 
course of the Edinburgh meeting on nonstandard analysis 
(August 1996). 
The authors are thankful to the organizers of the meeting for 
the opportunity to give a preliminary talk. 

%By {\it simple\/} (polygon, curve) we shall always mean: 
%non--self--intersecting.

\section*{Comment}

We shall use the $\IST$ language of nonstandard analysis 
(see \textsc{Nelson}~\cite{ne77}) to present the proof. However, 
as the reasoning involves only some very basic nonstandard 
notions, the exposition will be equally well understood by the 
followers of the ``asterisk'' version of nonstandard analysis 
(although nonstandard polygons should be called 
{\it hyperpolygons\/} or \dd\ast{\it polygons\/}, $\cK$ should 
be sometimes replaced by $\aK$ \etc)

\section*{Plan of the proof}

Starting the proof of the Jordan theorem, we consider a simple 
closed curve $\cK =\ans{K(t):0\leq t<1}$ where 
$K:\dR\,\lra\,\dR^2$ is a standard continuous \dd1 periodic 
function which is injective modulo $1$ (\ie\ $K(t)=K(t')$ implies 
$t-t'\equiv 0\bmod 1$). From $K(t)\approx K(t')$ it follows then 
that $t\approx t'\;\bmod1$.\vom

{\it Section \ref{appr}.} 
%Working in a fixed nonstandard domain, 
We infinitesimally approximate $\cK$ by a simple polygon $\Pi,$ 
using a construction, due to Jordan, of consecutive cutting 
loops in an originally self-intersecting approximation.\vom

{\it Section \ref{domain}.} We define the interior $\cK\inte$ 
as the open set of all (standard) points which belong to 
$\Pi\inte$ but does not belong to the monad of $\Pi.$  
(The Jordan theorem for polygons is taken for granted; this 
attaches definite meaning to $\Pi\inte$ and $\Pi\exte$.) 
$\cK\exte$ is defined accordingly.\vom

{\it Section \ref{nonempt}.} To prove that $\cK\inte$ is nonempty, 
we take a longest diameter $AB$ of $\Pi,$ draw two straight lines, 
$\al$ and $\ba,$ through resp.\ $A$ and $B$ orthogonally to $AB,$ 
and a parallel line $\ga$ between $\al$ and $\ba$ at equal distance 
from them. Now $\Pi$ is divided by $A$ and $B$ on two simple 
disjoint broken lines, $\cL$ (the left part) and $\cR$ (the right 
part). The straight segment $CD$ of $\ga$ bounded by the rightmost 
intersection $C$ of $\ga$ with $\cL$ and the next to the right 
intersection $D$ of $\ga$ with $\cR$ is included in $\Pi\inte$ 
and contains a point which does not belong to the 
monad~of~$\Pi$.\vom

{\it Section \ref{pconn}.} To prove that $\cK\inte$ is 
path-connected we define a simple polygon $\Pi'$ which lies 
entirely within $\Pi\inte,$  
does not intersect $\cK,$ and contains all points of $\cK\inte.$ 
This easily implies the path-connectedness.

\section{Approximation by a simple polygon}
\label{appr}

We say that an (internal) polygon $\Pi=P_1P_2\ldots P_nP_1$ 
($n$ may be infinitely large) 
\emph{approximates\/} $\cK$ if there is an internal sequence 
of reals $0\leq t_1<t_2<\ldots <t_n<1$ such that 
\ben
\itemsep=1mm
\item[$(\dag)$] $P_i=K(t_i)$ for $1\leq i<n,$ \hfill and \hfill $\,$

\item[$(\ddag)$] $t_n-t_1>\frac 12$ and 
$t_{i+1}-t_i<\frac 12$ for all $1\leq i<n$.
\een
%We define in this case 
%$\Da(\Pi)=\max_{1\leq k\leq n}\left| P_kP_{k+1}\right|$ 
%(where it is understood that $P_{n+1}=P_1$.)
%
We say that $\Pi$ \emph{approximates\/} $\cK$ 
\emph{infinitesimally\/} if in addition $\Da(\Pi)\approx 0,$ 
where $\Da(\Pi)=\max_{1\leq k\leq n}\left| P_kP_{k+1}\right|$ 
(it is understood that $P_{n+1}=P_1$).    

\ble
\label{app}
Let\/ $\Pi=P_1\ldots P_nP_1$ approximate\/ $\cK$ 
infinitesimally. Then
\ben
\itemsep=1mm
\def\theenumi{(\roman{enumi})}
\def\labelenumi{{\rm \theenumi}}
\item\label i  
$\msur n$ is infinitely large, $t_{i+1}\approx t_i$ for all 
$1\leq i<n,$ $t_1\approx 0,$ and $t_n\approx 1\,;$ 

\item\label {ii}  
%there exists an infinitesimal\/ $\vep>0$ such that\/ 
%$\Pi $
%is within the\/ \dd\vep nbhd of\/ $\cK$ and vice versa\/ $\cK$
%is within the\/ \dd\vep nbhd of\/ $\Pi$.
$\hsur\mon\cK=\mon\Pi$.

\item\label n
If\/ $k<l$ and\/ $P_k\approx P_l$ then$:$  
\underline{either}\/ $t_k\approx t_l$ and the arcs\/ $P_kP_l$ of 
both\/ $\cK$ and\/ $\Pi$ are contained in $\mon P_k=\mon P_l
\;;$
\underline{or}\/ $t_k\approx 0,\msur$ $t_l\approx 1,$ and the 
arcs\/ $P_lP_k$ 
of both\/ $\cK$ and\/ $\Pi$ are contained in $\mon P_k
=\mon P_l
$.
%\een
\een
\ele
\bpf 
\ref i\ 
The requirement $(\ddag)$ does not allow the reals $t_k$ to 
collapse into a sort of infinitesimal ``cluster'' or into a pair 
of them grouped around $0$ and $1;$  
%pair of $t_1=0$ and $t_2=1$ (for $n=2$); 
both the cases are compatible with the conjunction of 
$(\dag)$ and $\Da(\Pi)\approx 0$.\vom

\ref{ii}\ 
$\da_i=\max_{t_i\leq t\leq t_{i+1}}\left| K(t)-K(t_i)\right|$ is
infinitesimal for each $1\leq i\leq n$ and therefore 
$\vep=2\max_{1\leq i\leq n}\da_i$ is infinitesimal and proves 
the assertion. 
\epf

\ble
\label{apol}
There exists a simple polygon which infinitesimally approximates 
$\cK$.
\ele
\bpf 
Taking $t_i=\frac in$ for some infinitely large $n$ results
in a polygon which infinitesimally approximates $\cK$. But 
it may have self-intersections.

First of all two adjacent sides, $P_kP_{k+1}$ and $P_{k+1}P_{k+2}$ 
may have an ``illegal'' self-intersection other than the common 
vertex $P_{k+1},$ \ie\ either $P_k$ is an inner point of 
$P_{k+1}P_{k+2}$ or $P_{k+2}$ is an inner point of $P_kP_{k+1}.$ 
In this case we simply eliminate the vertex $P_{k+1},$ so that 
the polygon $P_1\dots P_kP_{k+1}P_{k+2} P_nP_1$ reduces to 
$P_1\dots P_kP_{k+2}\dots  P_nP_1$.

Assume two non-adjacent sides intersect, \ie\ $P_iP_{i+1}$ intersects 
$P_jP_{j+1}$ for some $1\leq i<j-1<n.$ By the triangle inequality 
the shorter of the segments $P_iP_j$ and $P_{i+1}P_{j+1}$ is not 
longer than the longer of the segments $P_iP_{i+1}$ and 
$P_jP_{j+1}$ which is bounded in length by $\Da(\Pi)$. 

If $\left| P_iP_j\right| \leq \left| P_{i+1}P_{j+1}\right|$ then 
we consider that one of the polygons 
\dm
P_1\ldots P_iP_j\ldots P_nP_1\hspace{5mm}\mbox{and}\hspace{5mm}
P_iP_{i+1}\ldots P_jP_i
\dm
which is parametrically longer, \ie\ the first one
if $t_j-t_i\<\frac 12$ or the second one if $t_j-t_i>\frac 12.$ 
In the case when  
$\left| P_{i+1}P_{j+1}\right| \leq \left|P_iP_j\right|$ 
we take that one of the polygons 
$P_1\ldots P_{i+1}P_{j+1}\ldots P_nP_1$ and 
$P_{i+1}P_{i+2}\ldots P_{j+1}P_{i+1}$ which is
parametrically longer. ($\dots$ means that all indices in 
between are involved.)

In all the cases the resulting polygon $\Pi\new$ still 
approximates $\cK$ (Lem\-ma~\ref{app}\ref i easily implies 
$(\ddag)$ for $\Pi\new$)  and satisfies $\Da(\Pi\new)\<\Da(\Pi)$ 
because the only new side is not longer than a certain side of 
$\Pi$. 

This (internal) procedure does not necessarily reduce the number of 
self-intersections because for the one which is removed there may 
be others appearing on the newly introduced side of the reduced 
polygon $\Pi\new.$ But the number of vertices of 
$\Pi\new$ is strictly less than that of $\Pi$. 

Therefore the internal series of polygons arising from 
$\Pi$ by iterated applications of this reduction procedure 
eventually ends by a simple polygon $\Pi'.$ This polygon 
approximates $\cK$ and satisfies $\Da(\Pi')\<\Da(\Pi)$ by the 
construction, so $\Da(\Pi')\approx 0.$ Finally, as the 
(internal) requirement $(\ddag)$ is preserved at each step of 
the procedure, we conclude, by internal induction, 
that $\Pi'$ still satisfies $(\ddag),$ hence approximates 
$\cK$ infinitesimally.\epf

\section{Definition of the inner and outer region}
\label{domain}

Let us fix for the remainder a polygon 
$\Pi =P_1P_2\ldots P_nP_1$ which approximates $\cK$ infinitesimally.

Let $I$ be the open set of all (standard) points $A\in\Pi\inte$ 
which have a non-infinitesimal distance from $\Pi.$ We put 
$\cK\inte=I$ and call this
the \emph{inner region} of the curve $\cK.$ In the same way we
define the open set $E$ of all (standard) points from 
$\Pi\exte$ which have non-infinitesimal distance from $\Pi,$ put 
$\cK \exte=E$ and call this the \emph{outer region} of the curve 
$\cK$.

Omitting rather elementary proofs that $\cK\inte$ is bounded,  
$\cK\exte$ is unbounded, and the complement of the union of both 
sets equals the curve $\cK,$ 
%Let $P$ be a (standard) point of that complement. 
%Then any standard circle around $P$ intersects $\Pi.$ By overspill 
%there is a point $P'\in \Pi$ with $P\approx P'.$ This implies 
%by Lemma \ref{app} that there is a point $P''\approx P$ on $\cK .$ 
%Then $P$ equals $P''$ as both are standard.  
%
let us prove that for $A\in \cK\inte$ and $B\in \cK\exte$ any 
(standard) continuous arc from $A$ to $B$ intersects $\cK.$ 
Indeed the arc must intersect $\Pi$ in some point $P$ 
because it starts in $\Pi\inte$ and ends in $\Pi\exte.$ (The Jordan 
theorem for polygons, transferred to the nonstandard domain, is 
applied.) By Lemma \ref{app}, there is a (standard) point 
$P'\in\cK$ infinitesimally close to $P\in\Pi.$ As $\cK$ and the 
arc are standard and closed, $P'$ is in $\cK$ and the arc.

\section{The inner region is non-empty}
\label{nonempt}

We prove that the inner region is not empty (that the 
outer region is not empty is trivial). 
%This is the hard part of the classical proof of \textsc{Jordan.}
%
By definition it suffices to prove the existence of a point in 
$\Pi\inte$ having non-infinitesimal distance from $\Pi$.

Arguing in the nonstandard domain, we let $A$ and $B$ be the two 
vertices of $\Pi$ with maximal distance between them; the distance 
by necessity must be non-infinitesimal. Let $\al$ and $\ba$ be two 
straight lines through $A$ and $B$ respectively and orthogonal 
to the segment $AB.$ The polygon $\Pi$ now consists of two 
broken lines $\cL$ and $\cR$ joining $A$ and $B,$ having no common 
points and not intersecting $\al$ and $\ba$ except for the points  
$A$ and $B$.

Let $\ga$ be a straight line parallel to $\al$ and $\ba$ and 
drawn between them at equal distance from both. Suppose that 
$\cL$ is the arc (among $\cL$ and $\cR$) first 
encountered when coming from left along $\ga.$ There must be a 
rightmost intersection point $C$ for $\ga$ and 
$\cL.$ Continuing from there further to the right there must be 
a first intersection $D$ with $\cR.$ One easily sees that the 
points between these two must be contained in the interior 
region $\Pi\inte$. 

Consider a point $E$ in the segment $CD$ which has the 
equal distance $d=\dist(E,\cL)-\dist(E,\cR)$ from both $\cL$ and 
$\cR.$ Note that $d$ is not infinitesimal. Indeed 
otherwise there are vertices $L\in \cL$ and $R\in \cR$ such that 
$L\approx E\approx R,$ 
%This would imply that either $E\approx L\approx A\approx R$ or 
%$E\approx L\approx B\approx R$ 
which is impossible by Lemma~\ref{app}\ref n as $\ga$ has 
non-infinitesimal distance from $\al$ and $\ba$.

\section{Path-connectedness}
\label{pconn}

Let $A$ and $B$ be two (standard) points in $\cK\inte$. We have 
to prove that there is a broken line joining $A$ with $B$ and 
not intersecting $\cK.$ This is based on the following lemma.

\ble
\label{rl}
There exists a simple polygon\/ $\Pi'$ lying entirely 
within\/ $\Pi\inte,$ containing no point of\/ $\aK$ in\/ 
$\Pi'\inte,$ and containing 
both\/ $A$ and\/ $B$ in\/ $\Pi'\inte$.
\ele
The lemma clearly implies the path-connectedness: indeed, $A$ 
can be connected to $B$ by a broken line which lies within 
$\Pi'\inte$ therefore does not intersect $\aK.$ By Transfer 
we get s standard broken line which connects $A$ and $B$ and 
does not intersect $\cK,$ as required.

\bpf{}of the lemma. Let an infinitesimal $\vep>0$ be defined as 
in the proof of Lemma~\ref{app}\ref{ii}, so that $\aK$ is 
included in the \dd\vep nbhd of $\Pi$.

Note that each side of $\Pi$ is infinitesimal by definition.   
For any side $PQ$ of $\Pi$ 
we draw a rectangle of the size $(|\wed PQ|+4\vep)\ti(4\vep)$ so 
that the side $PQ$ lies within the rectangle at equal distance 
$2\vep$ from each of the four sides of the rectangle. 

Let us say that a point $E$ is the {\it inner intersection\/} of 
two straight segments $\sg$ and $\sg'$ iff $E$ is an inner point of 
both $\sg$ and $\sg',$ and $\sg\cap\sg'=\ans{E}.$ 
For any point $C\in\Pi\inte$ which is either a vertex of some 
of the rectangles above, or an inner intersection of sides of two 
different rectangles in this family --- let us call points of 
this type {\it special points\/} --- let $\wed C{C'}$ be a  
shortest straight segment which connects $C$ with a point $C'$ 
on $\Pi\,;$ obviously $CC'$ is infinitesimal.   

The parts of the rectangles lying within $\Pi$ and the 
segments $\wed C{C'}$ for special points $C$ 
decompose the interior $\Pi\inte$ in a number of polygonal domains. 
Let $\Pi'$ be the polygon among them such that ${\Pi'}\inte$ 
contains $A.$ (Note that all the lines involved lie in 
the monad of $\Pi,$ hence none of them contains $A$ or $B$.) It 
remains to prove that ${\Pi'}\inte$ also contains $B,$ the other 
point.

It is clear that each vertex of $\Pi'$ is a special point (in the 
sense above). It follows that each side of $\Pi'$ is a part of 
either a side of one of the rectangles covering $\Pi$ or a segment 
of the form $\wed{C}{C'}$ --- therefore it is infinitesimal. 

Let $\Pi'=C_1C_2\dots C_n.$ We observe that, for any $k=1,...,n,$ 
the shortest segment $\sg_k=\wed{C_k}{C'_k},$ connecting $C_k$ with 
a point $C'_k$ in $\Pi,$ does not intersect ${\Pi'}\inte$ by 
construction. Moreover, by the triangle equality, the segments 
$\sg_k$ have no inner intersections. Therefore we may suppose that 
any two of them intersect each other only in such a manner that 
either the only intersection point is the common endpoint 
$C'_k=C'_l$ or one of them is an end-part of the other one.
Then the segments $\sg_k$ decompose the ring-like polygonal region 
$\cR$ between $\Pi$ and $\Pi'$ into 
$n$ open domains $\cD_k\msur$ $(k=1,...,n)$ defined as follows. 

If $\sg_k$ and $\sg_{k+1}$ are disjoint (${\sg}_{n+1}$ 
equals ${\sg}_1$) then the border of $\cD_k$ 
consists of $\sg_k,\msur$ $\sg_{k+1},$ the side $\wed{C_k}{C_{k+1}}$ 
of $\Pi',$ and that part $\curve{C'_k}{C'_{k+1}}$ of $\Pi$ 
%such that this polygon does not contain $A$ in its interior region. 
which does not contain any of $C'_l$ as an inner point. 
%If $C_k$ is an inner point of 
%$\sg_{k+1}$ and the segments have no other common points then 
%the border includes $\sg_k,$ the part $\wed{C_k}{C'_{k+1}}$ 
%of $\sg_{k+1},$ and $\curve{C'_k}{C'_{k+1}}.$ 
If $\sg_k$ and $\sg_{k+1}$ have the common endpoint $C_k=C_{k+1}$ 
and no more common points then the border shrinks to $\sg_k,\msur$ 
$\sg_{k+1},$ and $\wed{C_k}{C_{k+1}}.$ If, finally, one of 
the segments is an end-part of the other one then $\cD_k$ is 
empty. 

If now $B\in {\Pi'}\exte$ then $B$ belongs to one of 
the domains $\cD_k.$ If this is a domain of the first type 
then the infinitesimal simple arc $C'_kC_kC_{k+1}C'_{k+1}$  
separates $A$ from $B$ within $\Pi,$ which easily implies, 
by Lemma~\ref{app}\ref n, that either $A$ or $B$ belongs to 
$\mon \Pi,$ which is a contradiction with the choice of the 
points. If $\cD_k$ is a domain 
of second type then the barrier accordingly shrinks, with the 
same contradiction. 
%The domain of the third type is impossible:  
%$B$ would have distance ${\<\rho\sqrt 2}$ from $\Pi$).
\epf

\end{document}